\documentclass{amsart}
\newtheorem{theorem}{Theorem}[section]
\newtheorem{lemma}[theorem]{Lemma}
\newtheorem{proposition}[theorem]{Proposition}
\theoremstyle{definition}
\newtheorem{definition}[theorem]{Definition}

\theoremstyle{remark}
\newtheorem{remark}[theorem]{Remark}
\numberwithin{equation}{section}

\newcommand{\g}{{\overline{\gamma}}}
\DeclareMathOperator{\dist}{dist}

\DeclareMathOperator{\spa}{span}
\DeclareMathOperator{\ran}{ran}
\begin{document}

\title{Point evaluation and Hardy space
on a  homogeneous tree}

\author{Daniel Alpay
}

\address{Department of Mathematics,
Ben--Gurion University of the Negev, Beer-Sheva 84105, Israel}
\email{dany@math.bgu.ac.il}

\thanks{The first author
was supported by the Israel Science Foundation (Grant no.
322/00).}

\author{Dan Volok}
\address{Department of Mathematics,
Ben--Gurion University of the Negev, Beer-Sheva 84105, Israel}
\email{volok@math.bgu.ac.il}


\date{}


\begin{abstract}
We consider stationary multiscale systems as defined by
Basseville, Benveniste, Nikoukhah and Willsky. We show
that there are deep analogies with the discrete time non
stationary setting as developed by the first author, Dewilde and
Dym. Following these analogies we define a point evaluation with values
in a $C^*$--algebra and the corresponding ``Hardy space'' in
which Cauchy's formula holds. This point evaluation is used to
define in this context the counterpart of classical notions such
as Blaschke factors.
\end{abstract}

\maketitle
\section{Introduction}
The purpose of this paper is to explain relationships between the
theory of non-stationary linear systems indexed by the integers
and the theory of stationary linear systems indexed by a
homogeneous tree. We restrict ourselves to the case of input/ouput
systems and postpone the treatment of state space realizations to
a future publication. Our motivation for this study originates
with the works of Basseville, Benveniste, Nikoukhah and Willsky
(see \cite{MR95f:93057}, \cite{BNW}, \cite{MR95b:94008},
\cite{BBW}) where a general theory of multiscale system is
developed. In particular these authors distinguished special
families of stochastic processes (stationary and isotropic) and
developed  Schur-Levinson recursions for isotropic processes. They
also distinguish a special family of operators which they call
{\sl stationary transfer functions}. We associate to such
functions point evaluations whose domain and range are in a
$C^*$-algebra associated to the tree. We explain the relationships
with the point evaluation defined for upper-triangular operators
in \cite{MR93b:47027}. The analogies between the two cases allow
interfeeding. In particular, one can pose and solve interpolation
problems which should have applications to the theory of
multiscale systems. We plan to consider this in a future
publication. Some of the
results presented here were announced in \cite{a-volok1}.\\

The outline of the paper is as follows. In Section \ref{sec2} we briefly
review the setting developed in \cite{MR93b:47027}. In the third Section
we present part of the multiscale system theory relevant to the present
study. Although in part of a review nature, the point of view contains
some novelties, in particular by considering the $\ell_2$ space associated to
the tree. An important role is played by the so--called Cuntz relations (see
\cite{MR57:7189}). In the fourth Section we introduce the $C^*$-algebra of
constants and the Hardy space associated to the tree.
In Section \ref{sec5} we study the properties of the point evaluation and the notion of
Schur multipliers is studied in Section \ref{sec6}.\\

We would like to mention that similar techniques were used by Constantinescu
and Johnson (see
\cite{MR2003f:15039}, \cite{MR1929422}) in a different setting (tensor algebras
rather than Cuntz algebras).
Elements of system theory for Cuntz algebras
are being considered by Ball and Vinnikov in \cite{bv-c} and \cite{bv-new}.
Also, the connections between Cuntz algebras and wavelets are studied in
the book \cite{MR2003i:42001}.\\

{\bf Acknowledgments:} It is a pleasure to thank Prof. J. Ball for 
insightful comments on preliminary versions of this paper, and Prof. 
A. Benveniste who made available to us the reports \cite{BBW} and \cite{BNW},
which sparked our interest in the subject.

\section{The discrete time non stationary setting}
\label{sec2}
We briefly review the nonstationary setting developped in
\cite{MR92g:94002} and \cite{MR93b:47027}. We fix
a separable Hilbert space ${\mathcal N}$,
{\sl the coefficient space}, and consider
the Hilbert space  $\ell_{\mathcal N}^{2}$
of all two sided square summable sequences
$f= (f)_{i=-\infty}^{\infty}
  =(\ldots,f_{-1},\fbox{$f_{0}$},f_{1},\ldots)$ with components
$f_i \in {\mathcal N}$ provided with the standard inner product.
The set of bounded linear operators
from $\ell_{\mathcal N}^{2}$ into itself is denoted by ${\mathcal X}$.
Let $Z$ denote the bilateral backward shift operator
$$
\begin{array}{cccc}\left(Z f\right)_{i}
=f_{i+1}, & & i=\ldots,-1,0,1,\ldots . & \end{array}
$$
It is unitary on $\ell_{\mathcal N}^{2}$. An element $A\in
{\mathcal X}$ can be represented as an operator matrix $(A_{ij})$
with $A_{ij}=\pi^{\ast}Z^{i}AZ^{\ast j}\pi$ where $\pi$ denotes
the injection map: $u\in{\mathcal N}\mapsto
(\ldots,0,\fbox{$u$},0,\ldots)\in\ell_{\mathcal N}^{2}$. We denote
by ${\mathcal U}$ and ${\mathcal D}$ the spaces of
upper triangular and diagonal operators :
$$
{\mathcal U}=\left\{A\in{\mathcal X}: A_{ij}=0,\, i>j\right\},
\quad {\mathcal D}=\left\{A\in{\mathcal X}: A_{ij}=0,\, i\not
=j\right\}.
$$

For $W\in{\mathcal D}$ we set $W^{(j)}=Z^{\ast j}WZ^{j}$
and
\[
W^{[0]}=I, \quad W^{[n]}=WW^{(1)}W^{(2)}\ldots W^{(n-1)}=
(WZ^*)^nZ^n,\quad n\ge 1.
\]

\begin{theorem}\label{diagcan}
Let $F\in{\mathcal U}$ and $D\in{\mathcal D}$.
There exists
a unique sequence of operators $F_{[j]}\in{\mathcal D}$, $j=0,1,\ldots$ ,
namely $\left(F_{[j]}\right)_{ii}=F_{i-j,i}$, such that
\begin{equation}
\label{origdiag}
F=\sum_{n=0}^{\infty}Z^{n}F_{[n]}
\end{equation}
in the sense that
$F-\sum_{j=0}^{n-1}Z^{j}F_{[j]}\in Z^{n}{\mathcal U}$.
The operator
$\left(Z-W\right)^{-1}\left(F-D\right)$
belongs to ${\mathcal U}$ for $W\in{\mathcal D}$ such that
$r_{sp}\left(Z^{\ast}W\right)<1$ if and only if
\begin{equation}
D=\sum_{n=0}^{\infty}W^{[n]}F_{[n]}\stackrel{\rm def.}{=}F^\wedge(W).
\label{point-e}
\end{equation}
\end{theorem}

An operator $F=(F_{ij})\in {\mathcal X}$
 is a
Hilbert--Schmidt operator if all its entries $F_{ij}$ are Hilbert--Schmidt
operators on ${\mathcal N}$ and
$\sum_{ij}{\rm Tr}~F_{ij}^*F_{ij}<\infty$, where Tr stands for trace.
The set of these operators will be denoted by
${\mathcal X}_2$ and it is a Hilbert space with respect to the inner product
$$
\langle F,G\rangle_{{\mathcal X}_2}=\sum_{ij}{\rm Tr}~G_{ij}^*F_{ij}<\infty.
$$

The subspaces of upper triangular and diagonal operators which are moreover
Hilbert--Schmidt operators on ${\mathcal N}$
will be denoted by ${\mathcal U}_2$ and ${\mathcal D}_2$ .
The space ${\mathcal U}_2$ is a reproducing kernel Hilbert space with
reproducing kernel
$$\rho_W^{-1}=(I-ZW^*)^{-1}=\sum_0^\infty (ZW^*)^n=\sum_0^\infty
Z^nW^{[n]*}$$
in the sense that for all $W\in\Omega$, $E\in{\mathcal D}_2$,
and $F\in{\mathcal U}_2$, the operator $\rho_W^{-1}E\in{\mathcal U}_2$ and

\begin{equation}
\label{cauchy}
\langle F,\rho_W^{-1}E\rangle_{{\mathcal U}_2}={\rm Tr}~E^*F^\wedge(W).
\end{equation}

This last formula is the non stationary counterpart of Cauchy's formula for
Hardy functions.\\

The map $W\mapsto F^\wedge(W)$
(which was first introduced in \cite{MR92g:94002}) and its counterpart
when one replaces $(Z-W)^{-1}(F-D)$ by $(F-D)(Z-D)^{-1}$ in the above
theorem allow to translate most, if not all, of the classical analysis
of the Hardy space $H_2$ to the setting of upper triangular operators. The
analogue of the Hardy space $H_2$ is given by the Hilbert space of upper
triangular operators ${\mathcal U}_2$. See
\cite{AP1}, \cite{bgk-3}, \cite{DD-ot56} for sample applications.
As already mentioned and as we will see in the sequel, they have analogues in
the setting of multiscale system theory.

\section{Multiscale system theory}
Some of the fundamental notions in the classical theory of
discrete time linear systems are that of causality and
stationarity. In this section we    review the
 analogues of these notions, introduced
by Basseville, Benveniste, Nikoukhah and Willsky in the case of
multiscale linear systems.\\

Let \( \mathcal T\) be a homogeneous tree of order
\(q\geq 2\) - an infinite acyclic, undirected, connected graph
such that every node has exactly \(q+1\) branches (see
\cite{MR57:16426},
 \cite{MR54:12990}). We consider  a linear system of the form
\begin{equation}\label{mulsys}
g(t)=(Sf)(t),\end{equation}
where the input signal \(f=f(t)\) and the
output signal \(g=g(t)\) belong to the Hilbert space
\(\ell_2(\mathcal{T})\) of  square-summable sequences, indexed by
the nodes of \(\mathcal T\), and where \(S\)
is a bounded linear operator on \(\ell_2(\mathcal{T})\) (notation: 
\(S\in\mathbf X(\mathcal T)\)). Using the
notation \(\chi_t\) for the element of the standard basis of
\(\ell_2(\mathcal{T})\), supported at the node \(t\), one can
write
\begin{equation}\label{defioo}
(Sf)(t)=\sum_{u\in\mathcal T}s_{t,u}f(u),\text{ where }
s_{t,u}=[S\chi_u,\chi_t]\in\mathbb C\end{equation}
and where the sum is absolutely convergent by Cauchy--Schwarz inequality.\\

According to the standard tree compactification procedure, a
boundary point  of \(\mathcal T\) is an equivalence class
 of infinite paths
modulo finite number of edges.
 Let us distinguish a boundary point  of
\(\mathcal T\) and denote it  by
\(\infty_{\mathcal T}\). Then for each \(t\in\mathcal T\)
there exists a unique representative \({\mathcal T}^-_t\)
of the equivalence class
 \(\infty_{\mathcal T},\)
starting at \(t.\)  For a pair of nodes \(t,s\),
the paths \({\mathcal T}^-_t,{\mathcal T}^-_s\)
have to coincide after a finite number of edges;
 the first of their common nodes is denoted by \(s\wedge t\).
The  notion of   distance  \(\dist(s,t)\) between the nodes
\(s,t\), defined as the number of edges along the  path connecting
\(s\) and \(t\), allows to introduce the  partial order
\[
s\preceq t \quad\text{if}\quad \dist(s,s\wedge t)\leq
\dist(t,s\wedge t)
\]
and the  equivalence relation
\begin{equation*}
s\asymp t \quad\text{if}\quad \dist(s,s\wedge t)= \dist(t,s\wedge
t).
\end{equation*}
The equivalence classes, defined with respect to the equivalence
relation above, are called  horocycles.

\begin{definition}\label{defcausal}
The multiscale linear system \eqref{mulsys} is said to be causal
if for every node \(t\in\mathcal T\) the subspace
\[\{f\in\ell_2(\mathcal{T})\ :\  t\preceq \text{support}(f)\}\] is
\(S\)-invariant.
\end{definition}

In order to analyze Definition \ref{defcausal}, we consider the
primitive shifts on the tree. By convention they act on the
right and are defined as follows. The  primitive  upward shift
 \(\g:\mathcal{T}\mapsto\mathcal{T} \) is determined by \[
\forall t\in\mathcal{T}:\ t\g \preceq t, \ \dist(t\g ,t)=1.\] In
the choice of the primitive downward shifts there is  some
freedom; we assume that some such choice
\(\alpha_i:\mathcal{T}\mapsto\mathcal{T} ,\ 1\leq i\leq q,\)
\[\forall t\in\mathcal{T}:\ \{s\in\mathcal{T}\ :\ t\preceq s,\ \dist(t,s)=1\}=
\{t\alpha_1,\ldots,t\alpha_q\},\]
is fixed, as well. Furthermore we consider the primitive shift operators,
acting on the left on \(\ell_2(\mathcal{T})\)
and defined via convolution:
\[\g  f(t)=\frac{1}{\sqrt{q}}f(t\g ),\ \alpha_if(t)=f(t\alpha_i).\]
We compute the adjoint operators \(\gamma=\g^*\), \(\overline\alpha_i=
\alpha_i^*\):
\[\gamma f(t)=\frac{1}{\sqrt{q}}\sum_{s\g=t}f(s),\
\overline\alpha_if(t)=\left\{\begin{array}{l@{\quad }l}
 f(t\g),&t=t\g\alpha_i,\\
0,&\text{otherwise},
\end{array}
\right.
\]
and observe that the
following relations hold true:
\begin{align}\label{cuntz1}
\alpha_i\overline\alpha_j=\delta_{i,j},\quad &
\sum_{i=1}^q\overline\alpha_i \alpha_i=1,\\
\label{updec}\g=\frac{1}{\sqrt{q}}
\sum_{i=1}^q\overline\alpha_i,\quad & \gamma\g=1.\end{align}
Equations \eqref{cuntz1} are called the Cuntz relations.
Equation \eqref{updec} implies that
the primitive upward shift operator \(\g\) is an isometry from
\(\ell_2(\mathcal{T})\) into itself. However, it is not surjective
and thus \(\g\) is not unitary (see also \eqref{defdel} and
\eqref{delav} below). We also observe that
 for any pair of nodes \(t,u\in\mathcal T \)
there exist a unique choice of indices \(i_1,\ldots,i_n,j_1,\ldots,j_m\),
such that
\[t=(t\wedge u)\alpha_{ i_n}\cdots\alpha_{ i_1},\
 u=(t\wedge u)\alpha_{ j_m}\cdots\alpha_{ j_1}\]
(note that, according to the definition of \(t\wedge u\),
\(j_m\not=i_n\)).
Then for any \(f\in\ell_2(\mathcal{T})\) it holds that
\[f(u)=\overline\alpha_{i_1}\cdots\overline\alpha_{i_n}\alpha_{j_m}\cdots
\alpha_{j_1}f(t).\] This observation leads to a multiscale
analogue of Theorem \ref{diagcan}. It can be  formulated 
in terms of the point-wise convergence of a sequence of bounded operators on 
\(\ell_2(\mathcal T)\):
we shall say that a sequence of bounded operators \(S_n\) converges point-wise
to a bounded operator \(S\) if for every \(f\in\ell_2(\mathcal T)\)
and \(t\in\mathcal T\) \(\lim_{n\rightarrow\infty} (S_nf)(t)=(Sf)(t)\).
We note  that  on the dense subspace of
finitely supported functions  the point-wise convergence 
implies the convergence in strong operator topology.

\begin{theorem}\label{basic}
Any operator \(S\in\mathbf X(\mathcal T)\)
can be  represented as the point-wise converging series
\begin{equation}\label{ioo}
S=\sum_{n,m=0}^\infty\sum_{\substack{1\leq \substack{i_1,\ldots,i_n\\
 j_1,\ldots,j_m}\leq q\\j_m\not = i_n}}
\overline\alpha_{i_1}\cdots
\overline\alpha_{i_n}\alpha_{j_m}\cdots
\alpha_{j_1}S^{i_1,\ldots,i_n}_{j_1,\ldots,j_m},
\end{equation}
where  \(S^{i_1,\ldots,i_n}_{j_1,\ldots,j_m}\in \mathbf X(\mathcal T)\)
are diagonal operators, uniquely determined by
\begin{align}
\label{restrdiag}
&S^{i_1,\ldots,i_n}_{j_1,\ldots,j_m}\chi_t=0,\quad
t\not\in\mathcal{T}\alpha_{j_m}\cdots \alpha_{j_1},\\
\label{diagdef}
&S^{i_1,\ldots,i_n}_{j_1,\ldots,j_m}\chi_{t\alpha_{j_m}\cdots
\alpha_{j_1}}= [S\chi_{t\alpha_{j_m}\cdots
\alpha_{j_1}},\chi_{t\alpha_{i_n}\cdots\alpha_{i_1}}]
\chi_{t\alpha_{j_m}\cdots
\alpha_{j_1}},\quad  t\in\mathcal{T}.
\end{align}
\end{theorem}

\begin{proof}
Let \(f\in\ell_2(\mathcal{T})\) and \( t\in\mathcal{T}\) be fixed.
Using the operators, defined by \eqref{restrdiag}, \eqref{diagdef},
one can rewrite \eqref{defioo} as follows:
\begin{multline*}
Sf(t)=\sum_{u\in\mathcal T}[S\chi_u,\chi_t]f(u)\\
=\sum_{n,m=0}^\infty\sum_{\substack{1\leq
 j_1,\ldots,j_m\leq q\\t\g^n\alpha_{j_m}\not = t\g^{n-1}}}
[S\chi_{t\g^n\alpha_{j_m}\cdots \alpha_{j_1}},\chi_t]
f(t\g^n\alpha_{j_m}\cdots \alpha_{j_1})\\
=
\sum_{n,m=0}^\infty\sum_{\substack{1\leq \substack{i_1,\ldots,i_n\\
 j_1,\ldots,j_m}\leq q\\j_m\not = i_n}}
\overline\alpha_{i_1}\cdots\overline\alpha_{i_n}\alpha_{j_m}\cdots
\alpha_{j_1}S^{i_1,\ldots,i_n}_{j_1,\ldots,j_m}f(t),
\end{multline*}
and we obtain \eqref{ioo}, where the convergence is point-wise.\\

Furthermore, let \(t\in\mathcal T\) be fixed and let 
\(S\in\mathbf X(\mathcal T)\) be
 of the form \eqref{ioo}, where the
coefficients \(S^{i_1,\ldots,i_n}_{j_1,\ldots,j_m} \) are
diagonal, then
\begin{multline*}
S\chi_t=\sum_{n,m=0}^\infty\sum_{\substack{1\leq \substack{i_1,\ldots,i_n\\
 j_1,\ldots,j_m}\leq q\\j_m\not = i_n}}
\overline\alpha_{i_1}\cdots\overline\alpha_{i_n}\alpha_{j_m}\cdots
\alpha_{j_1}S^{i_1,\ldots,i_n}_{j_1,\ldots,j_m}\chi_t\\
=\sum_{n,m=0}^\infty\sum_{\substack{1\leq \substack{i_1,\ldots,i_n\\
 j_1,\ldots,j_m}\leq q\\j_m\not =
 i_n}}[S^{i_1,\ldots,i_n}_{j_1,\ldots,j_m}\chi_t,\chi_t]
\overline\alpha_{i_1}\cdots\overline\alpha_{i_n}\alpha_{j_m}\cdots
\alpha_{j_1}\chi_t\\
=\sum_{n,m=0}^\infty\sum_{\substack{1\leq
 i_1,\ldots,i_n\leq q\\i_n\not = \hat i_m}}
[S^{i_1,\ldots,i_n}_{\hat i_1,\ldots,\hat i_m}
 \chi_t,\chi_t]\chi_{t\g^m\alpha_{i_n}\cdots\alpha_{i_1}},
\end{multline*}
where the indices \(\hat i_k\) are determined by \(
t\g^k\alpha_{\hat i_k}=
t\g^{k-1}\).
Since the sum above is taken over \(i_n\not = \hat i_m\), all the
summands are mutually orthogonal and  \eqref{diagdef} follows.
Thus, under the restriction \eqref{restrdiag}, the coefficients
\(S^{i_1,\ldots,i_n}_{j_1,\ldots,j_m} \) are determined uniquely.

\end{proof}

\begin{proposition}\label{propcausal}
The multiscale linear system \eqref{mulsys} is causal
if, and only if, the coefficients of the
 representation \eqref{ioo} for \(S\) satisfy
\begin{equation}\label{charcausal}
S^{i_1,\ldots,i_n}_{j_1,\ldots,j_m}=0, \text{ whenever }
n<m.
\end{equation}
\end{proposition}

\begin{proof}
According to Definition \ref{defcausal}, the multiscale linear
system \eqref{mulsys} is  causal if and only if
\[[S\chi_u,\chi_t]=0,\text{ whenever }t\not\preceq u.\]
Hence \eqref{charcausal} follows immediately from \eqref{diagdef}.

\end{proof}

Next we turn to the notion of stationarity. As in the classical
case, this should mean translation-invariance. However, here
 the primitive downward shifts are one-to-one but not onto, while
the primitive upward shift is  onto, but not one-to-one. In particular,
neither is a tree isometry
(a tree isometry is a graph automorphism which preserves distances)
and hence is not suitable for the role of a translation. Instead, we shall
say that a tree isometry \(\tau:\mathcal T\mapsto \mathcal T\) is
 a  primitive translation  if for every \(t\in\mathcal T\)
\begin{equation}\label{trans}
t\tau\g\asymp t.
\end{equation}

Let us analyze  the structure of a primitive translation \(\tau\).
First of all, we note that, since \(\tau\) is a tree isometry,
\[\dist(t\g\tau,t\tau)=1.\] According to \eqref{trans},
\[t\g\tau\g\asymp t\g\preceq t\asymp t\tau\g\] and hence
\[t\g\tau\preceq t\tau.\] Therefore, by definition,
\(t\g\tau=t\tau\g\) and we conclude that \(\tau\) commutes with
the primitive upward shift:
\begin{equation}\label{comtrans}
\tau\g=\g\tau.
\end{equation}
Furthermore, we observe that \(\tau\g\) must have a fixed point.
Indeed, let \(t\in\mathcal T\) and let
\[v=t\wedge t\tau=t\g^n=t\tau\g^{n+1},\text{ where }
n=\dist(t\wedge t\tau,t).\] Then, by \eqref{comtrans},
\[v\tau\g=t\g^n\tau\g=t\tau\g^{n+1}=v.\]
Thus we obtain a unique sequence of nodes
\[v_j=v\tau^j,\quad j\in\mathbb Z,\]
satisfying
\[v_j\g=v_{j-1},\ v_j\tau=v_{j+1}.\]
It is called the skeleton of the primitive translation \(\tau\).
Each node \(v_j\) of the skeleton corresponds to the
non-homogeneous tree \({\mathcal T}_{v_j}^+\), which is the
maximal connected subgraph of \(\mathcal T\), satisfying
\[{\mathcal T}_{v_j}^+\bigcap\{v_{j-1},v_j,v_{j+1}\}=v_j.\]
It is mapped isometrically by \(\tau\) onto \({\mathcal
T}_{v_{j+1}}^+\).\\

The operator of convolution with a primitive translation is 
unitary on \(\ell_2(\mathcal T)\). By abuse of notation, we
denote both the group of the tree isometries, generated by
primitive translations, and the group of the corresponding
convolution operators by \({\bf A}(\mathcal T)\).
\begin{definition}\label{defstat}
The multiscale linear system \eqref{mulsys} is said to be
stationary if \(S\) commutes with every \(\tau\in{\bf A}(\mathcal
T)\).
\end{definition}

\begin{lemma}\label{lti}
The multiscale linear system \eqref{mulsys}  is stationary if and
only if the value of the scalar product
\([S\chi_u,\chi_t]\) depends only on  \(\dist(t\wedge u,t)\) and
\(\dist(t\wedge u,u)\).
\end{lemma}

\begin{proof}
First of all, we note that, according  to Definition
\ref{defstat}, the multiscale linear system \eqref{mulsys}  is
stationary if and only if for every pair of nodes
\(t,u\in\mathcal T\) and every \(\tau\in{\bf A}(\mathcal T)\)
\begin{equation}\label{interstat}
[S\chi_u,\chi_t]=[S\chi_{u\tau},\chi_{t\tau}].\end{equation}

So let us assume that   \([S\chi_u,\chi_t]\) depends only on
\(n=\dist(t\wedge u,t)\) and \(m=\dist(t\wedge u,u)\) and let
\(\tau\in{\bf A}(\mathcal T)\). Without loss of generality, we can
also assume that \(\tau\) is a primitive translation. Let us fix
now a pair \(t,u\), then it follows from \eqref{trans} that
\[(t\wedge u)\tau=t\g^n=u\g^m.\] Since
\[\dist(t\tau,u\tau)=\dist(t,u)=m+n,\] we
conclude that
\[(t\wedge u)\tau=(t\tau)\wedge(u\tau)\]
and, therefore, \eqref{interstat} holds true.\\

Conversely, let us assume that the multiscale linear system
\eqref{mulsys}  is stationary. Then, in view of \eqref{interstat},
it suffices to prove that for any two pairs of nodes \(t,u\) and
\(t^\prime,u^\prime\), satisfying
\[\dist(t\wedge u,t)=\dist(t^\prime\wedge u^\prime,t^\prime)=n,\
\dist(t\wedge u,u)=\dist(t^\prime\wedge u^\prime,u^\prime)=m,\]
there exists \(\tau\in{\bf A}(\mathcal T),\) such that
\[t\tau=t^\prime,\ u\tau=u^\prime.\]
Such an isometry can be constructed as follows. Denote
\[\dist((t\wedge u)\wedge(t^\prime\wedge u^\prime),t\wedge u)=k,\
\dist((t\wedge u)\wedge(t^\prime\wedge u^\prime),t^\prime\wedge
u^\prime)=p.\]
In the case \(t=u\) we choose arbitrary primitive
translations \(\tau_1\) and \(\tau_2\), such that \[ t\tau_1\g=t,\
t^\prime\tau_2\g=t^\prime,\]
 and set
\[\tau=\tau_1^{-k}\tau_2^{p}.\]
In the case
  \(t\not=u\) we assume, without loss of generality, that
\(n\not=0\) and choose \(\tau_1\) as above.
  Then
\[t\tau_1^{-n-k}=t\g^{-n-k}
=(t\wedge u)\wedge(t^\prime\wedge u^\prime)=t^\prime\g^{-n-p}.\]
Now let \(t^\prime\)  belong to the skeleton of a primitive
translation \(\tau_2\). Then
\begin{align*}
u\tau_1^{-n-k}\in{\mathcal T}_{t^\prime\tau_2^{-2n-p}}^+,&\quad
u^\prime\in{\mathcal T}_{t^\prime\tau_2^{-n}}^+,\\
\dist(u\tau_1^{-n-k},t^\prime\tau_2^{-2n-p})&=
\dist(u^\prime,t^\prime\tau_2^{-n})=m.\end{align*} Thus \(\tau_2\)
can be chosen so that
\[u\tau_1^{-n-k}\tau_2^{n+p}=u^\prime,\]
and we set
\[\tau=\tau_1^{-n-k}\tau_2^{n+p}.\]
\end{proof}

\begin{remark}\label{multo}
In view of the formula \eqref{diagdef}, Lemma \ref{lti} implies that 
the multiscale linear system \eqref{mulsys}  is stationary if, and
only if,  each coefficient \(S^{i_1,\ldots,i_n}_{j_1,\ldots,j_m}\)
in the  series \eqref{ioo} has a constant (except for normalizing
zeroes -- see \eqref{restrdiag}) diagonal and,
moreover, for \(t\in\mathcal{T}\alpha_{j_1}\cdots \alpha_{j_m}\)
the diagonal entry
\([S^{i_1,\ldots,i_n}_{j_1,\ldots,j_m}\chi_t,\chi_t]\) depends
only on \(n,m\). The first condition is the multiscale analogue
of the Toeplitz condition. Unlike the discrete time case, here it is 
weaker than the stationarity condition.
\end{remark}

\begin{theorem}\label{transinv}
The multiscale linear system \eqref{mulsys}  is stationary if, and
only if, 
\[S\in\overline\spa_{\mathbb C}\{\g^n\gamma^m:n,m\in\mathbb Z_+\},\]
where the closure is taken in the point-wise sense.
In this case, the  system \eqref{mulsys}  is also causal if, and
only if, 
\[S\in\overline\spa_{\mathbb C}\{\g^n\gamma^m:n\geq m\}.\]
\end{theorem}

\begin{proof} 
Let us assume first that  the
multiscale linear system \eqref{mulsys}  is stationary. Then,
since \[\alpha_{j_m}\cdots
\alpha_{j_1}=\alpha_{j_m}\cdots \alpha_{j_1}\pi_{j_1,\ldots,
j_m},\] where \(\pi_{j_1,\ldots, j_m}\) denotes the orthogonal
projection onto the subspace of signals supported in
\(\mathcal{T}\alpha_{j_m}\cdots \alpha_{j_1}\), Remark \ref{multo} means 
that the representation
\eqref{ioo} can  be rewritten  as
\[S=\sum_{n,m\in\mathbb Z^+}\sum_{\substack{1\leq \substack{i_1,\ldots,i_n\\
 j_1,\ldots,j_m}\leq q\\j_m\not = i_n}}
\overline\alpha_{i_1}\cdots\overline\alpha_{i_n}\alpha_{j_m}\cdots
\alpha_{j_1} s_{n,m},\quad s_{n,m}\in\mathbb C.\] 
But the partial sums of this  series belong to 
\(\spa_{\mathbb C}\{\g^n\gamma^m:n,m\in\mathbb Z_+\}\), since
 \eqref{updec} leads to
\[\sum_{\substack{1\leq \substack{i_1,\ldots,i_n\\
 j_1,\ldots,j_m}\leq q\\j_m\not = i_n}}
\overline\alpha_{i_1}\cdots\overline\alpha_{i_n}\alpha_{j_m}\cdots
\alpha_{j_1}=q^{\frac{m+n}{2}}\g^{n}
\gamma^m-q^{\frac{m+n-2}{2}}\g^{n-1}\gamma^{m-1}.\]
If the system \eqref{mulsys} is also causal then,
according to Proposition \ref{propcausal}, we have
\(s_{n,m}=0\) for \(n<m\), hence
\[S\in\overline\spa_{\mathbb C}\{\g^n\gamma^m:n\geq m\}.\]

In order to prove the converse statements, we note first
that if \(S\) is of the form
\(S=\g^n\gamma^m\) then, because of \eqref{comtrans}
and the fact that
\[\tau\in\mathbf A(\mathcal T)\implies
\tau^*=\tau^{-1}\in\mathbf A(\mathcal T),\]
the system \eqref{mulsys} is stationary. In the case
\(n\geq m\) it is also causal, as follows from
Proposition \ref{propcausal}. It only remains to
observe that, in view of our Definitions \ref{defcausal},
\ref{defstat} and the fact that each operator 
\(\tau\in\mathbf A(\mathcal T)\) is a convolution operator,
 the properties of causality and stationarity are preserved
when passing to the point-wise limit. 

\end{proof}

\section{Stationary multiscale systems and non-stationary
discrete time systems}
The main goal of the present work is to investigate the multiscale
systems of the form \eqref{mulsys} which are both causal and stationary.
We denote the Banach algebra of  corresponding operators $S$ by \({\bf
U}(\mathcal T)\). According to  Theorem
\ref{transinv}, 
\begin{equation}
\label{sigman}
\mathbf U(\mathcal T)=
\overline\spa_{\mathbb C}\{\g^n\sigma_m:n,m\in\mathbb Z_+\},
\end{equation}
where the closure is taken in the point-wise sense
and
 \begin{equation}\label{defdel}
\sigma_{m}=\g ^{m}\gamma^{m},\quad m\in\mathbb Z_+.
\end{equation}
We note that the self-adjoint operator \(\sigma_{m}\in{\bf U}(\mathcal T),\)
defined by \eqref{defdel}, is just the average
\begin{equation}\label{delav}
\sigma_mf(t)=\frac{1}{q^m}
\sum_{\substack{ u\asymp
t\\ \dist(u,t)\leq 2m }}f(u)\end{equation}
and hence 
acts on the subspace of functions, supported in a horocycle.
Furthermore,  since \(\g\) is an isometry,
\begin{equation}
\label{proddel} \sigma_m\sigma_n=\sigma_{m\vee n},
\end{equation}
 where \( m\vee
n\) denotes the maximal of two integers \(m\) and \(n\). Hence  it
is more convenient to consider 
\begin{equation}
\label{defom}\omega_m=\sigma_m-\sigma_{m+1}, \quad m\in\mathbb Z_+.
\end{equation}

\begin{theorem}\label{project}
The space \(\ell_2(\mathcal T)\) admits the orthogonal
decomposition
\begin{equation}\label{ordec} \ell_2(\mathcal
T)=\bigoplus_{m=0}^\infty{\mathcal W}_m, \end{equation}
where the orthogonal projections \(\omega_m:\ell_2(\mathcal
T)\mapsto{\mathcal W}_m\) are given by \eqref{defom}. Each subspace 
\({\mathcal W}_m\) is mapped  isometrically
onto \(\mathcal W_{m+1}\) by the upward shift operator \(\g\).
\end{theorem}

\begin{proof}
It follows immediately from \eqref{defdel}
\eqref{proddel} and \eqref{defom} that
\[\omega_m^*=\omega_m,\
\omega_m \omega_n=\delta_{m,n}\omega_n.\] 
In order to show that
\(\sum_{m=0}^\infty\omega_m\) converges in the strong operator topology
to the identity operator \(I\),
it suffices to  note that
\[\sum_{m=0}^n\omega_m =I-\sigma_{n+1}  \]
and that,   in view of \eqref{delav}, the sequence
\(\sigma_n\) converges in the strong operator topology
to \(0\).\\

Finally, we have to show that 
\[\g{\mathcal W}_{m}={\mathcal W}_{m+1}, \quad m\in\mathbb Z_+.\] 
But it follows from \eqref{updec},
\eqref{defdel} and \eqref{defom} that
\begin{equation}\label{commg}
\omega_{m+1}=\g\omega_m\gamma \text{ and }
\g\omega_m=\omega_{m+1}\g.\end{equation}
\end{proof}

It follows from Theorem \ref{project}
that an operator \(S\in\mathbf X(\mathcal T)\)
can be viewed as an operator matrix with blocks
corresponding to the orthogonal decomposition
\eqref{ordec} and thus the corresponding multiscale linear
system \eqref{mulsys} 
can be treated as a (non-stationary, in general)
discrete time system (see  \cite{MR93b:47027}
and Section \ref{sec2}). Let us investigate the causal stationary 
multiscale linear
systems from this point of view.

\begin{proposition}
\label{newprop1}
Let \(S\in\mathbf X(\mathcal T)\). Then \(S\in\mathbf U(\mathcal T)\)
if, and only if,
\begin{equation}
\label{oblock}
\omega_m S \omega_n=\left\{\begin{array}{l@{\quad}l}\g^{m-n}\omega_n
s_{m,n},&
m\geq n,\\
0,&\text{otherwise},\end{array}\right.
\end{equation}
where \(s_{m,n}\in\mathbb C\). 
\end{proposition}

\begin{proof}
Let us assume that \(S\in\mathbf U(\mathcal T)\). Then,
by \eqref{sigman}, there exists
a sequence \(S_k\in\spa_{\mathbb C}\{\g^j\sigma_i:i,j\in\mathbb Z_+\}\)
which converges to \(S\) point-wise. \eqref{commg}
implies that 
\[\omega_m S_k \omega_n=\left\{\begin{array}{l@{\quad}l}\g^{m-n}\omega_n
s_{m,n,k},&
m\geq n,\\
0,&\text{otherwise},\end{array}\right.\]
where \(s_{m,n,k}\in\mathbb C\). Since for every \(m\in\mathbb Z_+\)
and \(t\in\mathcal T\) \(\omega_m\chi_t\) has a finite support,
\[\lim_{k\rightarrow\infty}\omega_m S_k \omega_n=\omega_m S \omega_n
\text{ point-wise.}\] 
In particular, there exists  \(\lim_{k\rightarrow\infty}s_{m,n,k},\)
which we can denote by \(s_{m,n}\) and thus  obtain
\eqref{oblock}.\\

Conversely, if \eqref{oblock} holds, then
by Theorem \ref{project} we can represent \(S\) as the 
strongly converging series
\[S=\sum_{m=0}^\infty\omega_m S=
\sum_{m=0}^\infty\sum_{n=0}^m\omega_m S\omega_n.\]
Since the convergence in the strong operator topology
 implies the point-wise convergence, it suffices
to apply Theorem \ref{transinv} to complete the proof.

\end{proof}

From Proposition \ref{newprop1}  it follows
that a multiscale system \eqref{mulsys} is stationary and causal
if, and only if,
\(S\) is "triangular"
with respect to the orthogonal decomposition \eqref{ordec} -- that is,
for every \(n\in\mathbb Z_+\) the subspace of piece-wise constant functions
\(\oplus_{m\geq n}\mathcal{W}_m=\ran\sigma_n\)
is \(S\)-invariant\footnote{In the language of 
nest algebras (see \cite{MR52:3979} \cite{MR84e:93003})
$S$ is in the nest algebra associated with the nest $\{\oplus_{m\ge n}
{\mathcal W}_m:n\in\mathbb Z_+\}$.} -- and, moreover, the blocks of \(S\) 
are complex constants. However, the subspaces \(\mathcal W_m\)
are infinite-dimensional and hence
 the only Hilbert-Schmidt element of \(\mathbf U(\mathcal T)\)
is \(0\). Nevertheless, we can adapt the techniques, developed
for the non-stationary discrete time systems, to the present setting.\\

In order to formulate the appropriate analogue of Theorem
\ref{diagcan}, we consider the space of operators
\[
\mathbb K=\{{\bf c}=\sum_{m=0}^\infty c_m\omega_{m}\ :\ c_m\in\mathbb C,
\sup_{m\in\mathbb Z_+}|c_m|<\infty,\}\]
where the convergence is in the strong operator topology.
According to Theorem \ref{project} and Proposition \ref{newprop1},
 \(\mathbb K\)
 is a subalgebra of \({\bf U}(\mathcal T)\)
and, moreover,  
a commutative \(\mathbb{C}^*\)-algebra\footnote{For
 background on \(\mathbb{C}^*\)-algebras we refer the reader to
\cite{MR39:7442}.}, isometric to
\(\ell_\infty(\mathbb
Z^+)\). For \(\mathbf c\in\mathbb K\) we shall
use the notation
\[\overline{\mathbf c}=\mathbf c^*=\sum_{m=0}^\infty \overline c_m\omega_{m}.\]
We also introduce the notion of the row-wise
(with respect to the orthogonal decomposition \eqref{ordec})
convergence: we shall say that
 a sequence of  \(S_n\in\mathbf X(\mathcal T)\) converges row-wise
to  \(S\in\mathbf X(\mathcal T)\) if for every \(m\in\mathbb Z_+\)
the sequence \(\omega_m S_n\) converges to \(\omega_m S\)
in the operator norm.

\begin{theorem}\label{ogo}
Let \(S\in\mathbf X(\mathcal T)\).
Then \(S\in{\bf U}(\mathcal T)\) if, and only if, it can be  represented as
a row-wise converging  series
\begin{equation}\label{formal1p}
S=\sum_{k=0}^{\infty}\g ^k{\bf s}_k,\quad \mathbf s_k\in\mathbb K.
\end{equation}
In this case
the operators \(\mathbf s_k\) are determined uniquely
by
\begin{equation}\label{defsk}
\omega_n\mathbf s_k=\gamma^k\omega_{n+k}S\omega_n,\quad n,k\in\mathbb
Z_+,
\end{equation}
and
it
holds that
\begin{equation}\label{funrel}\sum_{k=0}^\infty \overline{\bf s}_k{\bf
s}_k\omega_n=
\frac{q^{n+1}}{q-1}\|S\omega_n\chi_t\|^2\omega_n,\quad 
t\in\mathcal T,\ n\in\mathbb{Z}^+.
\end{equation}
\end{theorem}

\begin{proof}
First, let us assume that there exists
a sequence \(\mathbf s_k\in\mathbb K\)
such that the series \(\sum_k\g^k\mathbf s_k\) converges row-wise
to an operator \(S\in\mathbf X(\mathcal T)\).
Then for any \(m,n\in\mathbb Z_+\) 
\[\sum_{k=0}^\infty\omega_m\g^k\mathbf s_k\omega_n=\omega_m S\omega_n.\]
Then, since
\[\omega_m\g^k\mathbf s_k\omega_n=\delta_{k,m-n}\g^k\mathbf s_k\omega_n,\]
\(S\in\mathbf U(\mathcal T)\) by Proposition \ref{newprop1}. Moreover,
\[\g^k\mathbf s_k\omega_n=\omega_{n+k} S\omega_n,\quad n,k\in\mathbb
Z_+,\] hence \eqref{defsk} holds true.\\

Conversely, let us assume that \(S\in\mathbf U(\mathcal T)\)
and let \(m\in\mathbb Z_+\).
Then, according to Proposition \ref{newprop1}, 
we can define \(\mathbf s_k\in\mathbb K\) by \eqref{defsk}
and observe that
\[\omega_m S=\sum_{n=0}^m\omega_m S\omega_n=\sum_{n=0}^m\omega_m
\g^{m-n}\mathbf s_{m-n}=\sum_{n=0}^m\omega_m
\g^{n}\mathbf s_{n}=\sum_{n=0}^\infty\omega_m
\g^{n}\mathbf s_{n}.\]
Finally,
\begin{multline*}
\|S\omega_n\chi_t\|^2\omega_n 
=\sum_{m= n}^\infty\|\omega_mS\omega_n\chi_t\|^2\omega_n
=\sum_{m=n}^\infty\|\g^{m-n}\mathbf s_{m-n}\omega_n\chi_t\|^2\omega_n\\
=\sum_{m=0}^{\infty}\|{\bf s}_m\omega_n\chi_t\|^2\omega_n
=\|\omega_n\chi_t\|^2\sum_{k=0}^\infty\overline{\bf s}_k{\bf s}_k\omega_n
=\left(\dfrac{1}{q^n}-\dfrac{1}{q^{n+1}}\right)\sum_{k=0}^\infty
\overline{\bf s}_k{\bf s}_k\omega_n,
\end{multline*}
and we obtain \eqref{funrel}.

\end{proof}

Following the analogy with the non--stationary setting, we consider the
following ideal of
\(\mathbb K\):
\[
\mathbb K_2=\{{\bf c}\in\mathbb K\ :\
\sum_{m=0}^\infty|c_m|^2<\infty\}.\] It is
a Hilbert space, isometric to \(\ell_2(\mathbb{Z}^+)\).
We also consider
 the  \(\mathbb K\)-module\footnote{For background
on modules over a \(\mathbb C^*\)-algebra
see \cite{MR15:327f}, \cite{MR33:584},
\cite{wp}.}:
\[
{\bf H}_2(\mathcal T)=\{S=\sum_{k=0}^\infty\g^k{\bf s}_k\ :\
{\bf s}_k\in\mathbb K_2, \sum_{k=0}^\infty\|{\bf s}_k\|_2^2
<\infty \}.\]
At this point we consider the power series in the definition
above as formal. However, we shall see later on
(Proposition \ref{hsio})
that such a series converges in the operator norm. 
This is the analogue of the space of Hilbert-Schmidt operators
in the present setting.

\begin{proposition}
The \(\mathbb K\)-module
\({\bf H}_2(\mathcal T)\), considered as a vector
space over \(\mathbb C\) with the scalar product
\begin{equation}\label{inprod}
[F,G] =\sum_{k=0}^\infty[\mathbf f_k,\mathbf g_k],\end{equation}
 is a Hilbert space.
\end{proposition}

\begin{proof}
The proof is the same as in the case of Hilbert-Schmidt operators,
 hence we shall
give only an outline. Since the Cauchy--Schwarz
inequality  holds in \(\mathbb{K}_2\), the inner product
 \eqref{inprod} is well-defined in the whole of \({\bf H}_2(\mathcal T),\)
which is, therefore, a pre-Hilbert space. Hence, the Cauchy--Schwarz
inequality holds in \({\bf H}_2(\mathcal T)\)
as well, and
\begin{equation*}
\|F\|_2=\sqrt{[F,F]}\end{equation*}
is a well-defined norm. The completeness of \({\bf H}_2(\mathcal T)\)
with respect to this norm
can now be proved, using the triangle inequality.

\end{proof}

\begin{proposition}\label{hsio}
The following hold:
\begin{enumerate}
\item
The Hilbert space \({\bf H}_2(\mathcal T)\)
is contractively contained in \({\bf U}(\mathcal
T)\).
\item
Let \(S\in{\bf U}(\mathcal T)\) and \( t\in\mathcal T\). Then
\(S\in{\bf H}_2(\mathcal T)\) if, and only if,
\[\sum_{k=0}^\infty q^k\|S\omega_k\chi_t\|^2<\infty.\]
In this case the expression above is equal to
\((1-\frac{1}{q})\|S\|_2^2\).
\end{enumerate}
\end{proposition}

\begin{proof}
\mbox{}\\
\begin{enumerate}
\item
In view of completeness of \({\bf H}_2(\mathcal T),\) 
\(\mathbf U(\mathcal T)\)
 and
\(\ell_2(\mathcal T)\), it suffices to consider 
\(S\in\mathbf H_2(\mathcal T)\)  such that the coefficients \({\bf
s}_k\in\mathbb{K}_2\) are different from zero only for a finite
number of indices \(k\). Then \(S\in\mathbf U(\mathcal T)\)
and we have
\begin{multline*}
\qquad\quad\|Sf\|^2=\sum_{n=0}^\infty\|\omega_nSf\|^2
\leq\sum_{n=0}^\infty
\left(\sum_{k=0}^n\|{\bf
s}_k\omega_{n-k}f\|\right)^2\\
\qquad\leq\sum_{n=0}^\infty\left(\sum_{k=0}^n
\|\mathbf s_{k}\omega_{n-k}\|^2\right)
\left(\sum_{m=0}^n\|\omega_{n-m}f\|^2\right)\\
\leq\sum_{m=0}^\infty\|\omega_mf\|^2\sum_{n,k=0}^\infty
\|\mathbf s_{k}\omega_n\|^2=\|f\|^2\|S\|_2^2.
\end{multline*}
\item It follows from \eqref{funrel} that
\[\|S\|_2^2=\sum_{k=0}^\infty\|{\bf
s}_{k}\|_2^2=\frac{q}{q-1}\sum_{k=0}^\infty q^k\|S\omega_k\chi_t\|^2,\]
whenever either right-hand side or left-hand side is finite.
\end{enumerate}

\end{proof}

\begin{remark}
A consequence of Proposition \ref{hsio} is that
 \({\bf
H}_2(\mathcal T)\) is a left ideal in \({\bf U}(\mathcal T)\)
and, moreover, for any \(S\in{\bf U}(\mathcal T)\) and \(F\in{\bf
H}_2(\mathcal T)\) the inequality
\begin{equation}\label{multin}\|SF\|_2\leq\|S\|\|F\|_2\end{equation}
holds true.
\end{remark}

\section{Point evaluation}
\label{sec5}
In this  section we  exploit the analogy with the non-stationary
setting to associate the elements of \({\bf U}(\mathcal T)\)
with maps from \(\mathbb K\) into itself.\\

Recalling the identity \eqref{commg} and observing that
\(\omega_0\g=0\), we conclude that for every \({\bf
c}\in\mathbb{K}\) there exists an element
\begin{equation*}
{\bf c}^{(1)}=\sum_{n=0}^\infty \omega_n c_{n+1}\in{\mathbb K}
\end{equation*}
such that
\begin{equation*}
{\bf c}\overline{\gamma}=\g {\bf c}^{(1)}.
\end{equation*}
We note that
\(
\|{\bf c}^{(1)})\|\le\|{\bf c} \|\) and  \(({\bf
c}{\bf d})^{(1)}= {\bf c}^{(1)}{\bf d}^{(1)}\).
Furthermore,
we introduce the following notation:
\begin{align*}
{\bf c}^{(0)}={\bf c},&\quad
{\bf c}^{(n+1)}=\left ({\bf c}^{(n)}\right)^{(1)},\\
{\bf c}^{[0]}=1,&\quad
 {\bf c}^{[n+1]}={\bf c}^{[n]}{\bf c}^{(n)},\\
\rho({\bf c})&= \limsup_{n\rightarrow\infty} \|{\bf
c}^{[n]}\|^{\frac{1}{n}},\\
{\mathbb D}({\mathcal T})&=\{{\bf c}\in
{\mathbb K}\ :\ \rho({\bf c})<1\}.
 \end{align*}
The set  ${\mathbb D}({\mathcal T})$ is the counterpart of the open unit disk
in the present setting.

\begin{definition}\label{pointev}
Let \(S\in{\bf U}({\mathcal T})\) be given.
 For \({\bf c}\in\mathbb
{D}({\mathcal T})\) we define the point evaluation of \(S\) at
\({\bf c}\) by
\begin{equation}\label{defpoint}
S({\bf c})=\sum_{n=0}^\infty {\bf c}^{[n]}{\bf s}_n.
\end{equation}
\end{definition}

We note that \eqref{defpoint} is the analogue of \eqref{point-e}. 
We claim that the point evaluation \eqref{defpoint} is
well-defined. Indeed, the convergence
of the series \eqref{defpoint} in \(\mathbb K\) follows from
Theorem \ref{ogo}.
Also, if \(S\in{\bf U}({\mathcal T})\) is such
that for every \({\bf c}\in\mathbb {D}({\mathcal T})\) \(S({\bf
c})=0\), then, in particular,
\[S(0)=S(\omega_k)=S(\omega_k+\omega_{k+1})=\ldots=0.\]
Hence
\[{\bf s}_0=0,\ {\bf s}_1\omega_k=0,\ {\bf s}_2\omega_k=0,\ldots\]
and  \({\bf s}_n=0\) for \(n=0,1,2,\ldots\)
Furthermore,  if \(F\in{\bf H}_2({\mathcal T})\)
then for every \( {\bf c}\in\mathbb{D}(\mathcal T)\)
\(F({\bf c})\in\mathbb{K}_2\).
We list several other properties of the point evaluation
in the following

\begin{lemma}\label{proppoint}
The following hold:
\begin{enumerate}
\item
Let \(F,G\in{\bf U}({\mathcal T})\), \({\bf p}, {\bf q}\in\mathbb K\),
 \( {\bf c}\in\mathbb{D}(\mathcal T)\)
and assume
that \({\bf k}\in\mathbb K\) is invertible. Then
\begin{align}
\label{ppe0}(F{\bf p}+G{\bf q})({\bf c})&=F({\bf c}){\bf p}+
G({\bf c}){\bf q},\\
\label{ppe1}(FG)({\bf c})&=(F({\bf c})G)({\bf c}),\\
\label{ppe2}(\g^nF)({\bf c})&=
{\bf c}^{[n]}F({\bf c}^{(n)}),\\
\label{ppe3}({\bf k}F)({\bf c})&=F({\bf k}^{(1)}{\bf
k}^{-1}{\bf c}){\bf k}.
\end{align}
\item Let \(F\in{\bf H}_2({\mathcal T})\),
\( {\bf c}\in\mathbb{D}(\mathcal T)\). Then
\begin{equation}\label{bezou}
F-F({\bf c})=(\g-{\bf c})G,
\end{equation} where \(G\in{\bf H}_2({\mathcal T})\)
is given by
 \begin{equation}\label{21-decembre-2002} G=
\sum_{n,k=0}^{\infty}\g^n\left({\bf
c}^{(n+1)}\right)^{[k]}{\bf f}_{n+k+1}.
 \end{equation}
\end{enumerate}
\end{lemma}

\begin{proof} \ \\
\begin{enumerate}
\item
The relation \eqref{ppe0} follows immediately from Definition
\ref{pointev}. Furthermore,
\begin{multline*}
\qquad\quad(\g^nF)({\bf c})
=\left(\sum_{k=0}^\infty\g^{n+k}{\bf f}_k\right)({\bf c})=
\sum_{k=0}^\infty{\bf c}^{[n+k]}{\bf f}_k\\
\qquad=\sum_{k=0}^\infty{\bf c}^{[k]}\left({\bf c}^{[n]}\right)^{(k)}{\bf f}_k=
\left(\sum_{k=0}^\infty\g^k\left({\bf c}^{[n]}\right)^{(k)}{\bf f}_k\right)
({\bf c})\\
=\left({\bf c}^{[n]}\sum_{k=0}^\infty\g^k{\bf f}_k\right)
({\bf c})=(\g^n({\bf c})F)({\bf c}),\end{multline*}
and, in view of \eqref{ppe0}, we obtain \eqref{ppe1}.
Analogously,
\[
(\g^nF)({\bf c})=\sum_{k=0}^\infty{\bf c}^{[n+k]}{\bf f}_k
=\sum_{k=0}^\infty\left({\bf c}^{(n)}\right)^{[k]}{\bf f}_k{\bf c}^{[n]},
\]
and we obtain \eqref{ppe2}.
Finally, we note that
\[({\bf k}\g^n)({\bf c})=
{\bf c}^{[n]}{\bf k}^{(n)}=\left( {\bf c}{\bf k}^{(1)}{\bf k}^{-1}
\right)^{[n]}{\bf k},
\]
and  (\ref{ppe3}) follows from \eqref{ppe0}.
\item
First, we have to check that the series \eqref{21-decembre-2002}
belongs to \({\bf H}_2({\mathcal T})\).
But
\begin{multline*}
\qquad\quad\sum_{n=0}^{\infty}\|\sum_{k=0}^\infty\left({\bf
c}^{(n+1)}\right)^{[k]}{\bf f}_{n+k+1}\|_2^2\leq
\sum_{n=0}^{\infty}\left(\sum_{k=0}^\infty\|{\bf
c}^{[k]}\|\|{\bf f}_{n+k+1}\|_2\right)^2\\
\qquad\leq\sum_{n=0}^{\infty}\left(\sum_{k=0}^\infty\|{\bf
c}^{[k]}\|\right)\left(\sum_{m=0}^\infty\|{\bf
c}^{[m]}\|\|{\bf f}_{n+m+1}\|_2^2\right)\\
\leq\left(\sum_{k=0}^\infty\|{\bf
c}^{[k]}\|\right)^2\|F\|_2^2 <\infty,
\end{multline*}
 since \(\rho({\bf c})<1\). Hence, indeed,
\(G\in{\bf H}_2({\mathcal T})\).\\

Now we shall prove \eqref{bezou}. Without loss of generality,
we asssume that \(F=\g^m{\bf f}\), where \({\bf f}\in\mathbb{K}_2\).
 Then
\eqref{21-decembre-2002} means
\[G=\sum_{n=1}^{m}\g^{n-1}\left({\bf
c}^{(n)}\right)^{[m-n]}{\bf f}.\]
In particular, for \(m=1\)
we have \(G={\bf f}\) and \eqref{bezou} holds. For \(m\geq 2\)
we have
\[(\g^m-{\bf c}^{[m]}){\bf f}=(\g-{\bf c})\g^{m-1}{\bf f}
+(\g^{m-1}-{\bf c}^{[m-1]}){\bf f}
{\bf c}^{(m-1)}.\]
Applying induction on \(m\), we obtain
\begin{multline*}
\qquad\quad(\g^m-{\bf c}^{[m]}){\bf f}=(\g-{\bf c})\left(\g^{m-1}+
\sum_{n=1}^{m-1}\g^{n-1}\left({\bf
c}^{(n)}\right)^{[m-n-1]}{\bf c}^{(m-1)}\right){\bf f}\\
=(\g-{\bf c})\left(\g^{m-1}+\sum_{n=1}^{m-1}\g^{n-1}\left({\bf
c}^{(n)}\right)^{[m-n]}\right){\bf f}=(\g-{\bf c})G.
\end{multline*}
\end{enumerate}
\end{proof}

A consequence of Lemma \ref{proppoint} is that
 \(F\in{\bf H}_2({\mathcal T})\) satisfies \(F({\bf c})=0\)
if and only if
\(F\) is of the form \eqref{bezou}, where
\(G\in{\bf H}_2({\mathcal T})\).\\

Finally, we present an analogue of Cauchy's formula
(and of formula \eqref{cauchy}) for the space \({\bf H}_2(\mathcal T)\).
\begin{theorem}\label{RK}
Let \(F \in{\bf H}_2(\mathcal T)\), \({\bf c}\in\mathbb{D}(\mathcal T)\).
 Then for every \({\bf k}\in\mathbb{K}_2\) it holds that
\begin{equation}\label{repker}
[F({\bf c}),{\bf k}]=[F, K_{\bf c}{\bf k}],
\end{equation}
where \(K_{\bf c}\in{\bf U}(\mathcal T)\)
is given by
\begin{equation}\label{defkc}
K_{\bf c}=\sum_{n=0}^\infty\g^n\overline{{\bf c}}^{[n]}= (1-\g\
\overline{\bf c})^{-1}.
\end{equation}
\end{theorem}
\begin{proof}
Since \({\bf c}\in\mathbb{D}(\mathcal T)\), there exists
\(\epsilon\in(0,1)\), such that for \(n\) sufficiently large \(
\|{\bf c}\|\leq\epsilon^n\). It follows that the series
\eqref{defkc}  converges absolutely in  \({\bf U}(\mathcal T)\)
and  defines an
operator \(K_{\bf c}\), which satisfies
\[K_{\bf c}(1-\g\
\overline{\bf c})=(1-\g\ \overline{\bf c})K_{\bf c}=1.\] The
formula \eqref{repker} follows immediately from \eqref{defpoint}.

\end{proof}

\section{Schur multipliers}
\label{sec6}
Let us recall that, in view of Proposition \ref{hsio} and the subsequent
remarks, for any \(S\in{\bf U}(\mathcal T)\) the multiplication operator
\(\mathcal{M}_SF=SF\) is a bounded linear operator
from \({\bf H}_2(\mathcal T)\) into itself.

\begin{definition}
 \(S\in{\bf U}(\mathcal T)\)  is called a Schur multiplier if the
multiplication operator  \(\mathcal{M}_S\) is a contraction
in \({\bf H}_2(\mathcal T)\).
\end{definition}

\begin{theorem}
An element $S\in{\bf U}(\mathcal T)$ is
a Schur multiplier if and only if the
map ${\mathcal K}_S: {\mathbb D}({\mathcal T})\times{\mathbb D}({\mathcal
T})\mapsto {\mathbb K}$, defined by
\begin{equation}\label{defks}
\mathcal{K}_S({\bf c},{\bf d})=\sum_{n=0}^\infty {\bf c}^{[n]}\left(1-S({\bf
c})\overline{S({\bf d})}\right)^{(n)} {\overline{\bf d}}^{[n]}
\end{equation}
is positive:
for any\  \(m\in\mathbb{Z}^+\),
\({\bf c}_0,\ldots,{\bf c}_m\in\mathbb{D}(\mathcal T)\),
\({\bf k}_0,\ldots,{\bf k}_m\in\mathbb{K}_2\) it holds that
\begin{equation}\label{poscon}
\sum_{i,j=0}^m[\mathcal{K}_S({\bf c_i},{\bf c_j}){\bf k}_j,{\bf k}_i]\geq 0.
\end{equation}
\end{theorem}

\begin{proof}
Let \(S\in{\bf U}({\mathcal
T})\), \({\bf c}\in{\mathbb D}({\mathcal
T})\), \({\bf k}\in\mathbb{K}_2\). Then
\begin{multline*}\mathcal{M}_S^*(K_{\bf c}{\bf k})
=\sum_{n,m=0}^\infty\g^n\omega_m
[\mathcal{M}_S^*(K_{\bf c}{\bf k}),\g^n\omega_m]=
\sum_{n,m=0}^\infty\g^n
[K_{\bf c}{\bf k},S\g^n\omega_m]\omega_m\\=\sum_{n,m=0}^\infty
\g^n[{\bf k},(S\g^n\omega_m)({\bf
c})]\omega_m
=\sum_{n,m=0}^\infty
\g^n[{\overline{S({\bf
c})}}^{(n)}{\overline{\bf c}}^{[n]}{\bf k},\omega_m]\omega_m\\
=\sum_{n=0}^{\infty}\g^n{\overline{S({\bf
c})}}^{(n)}{\overline{\bf c}}^{[n]}{\bf k}.\end{multline*}
It follows that
\begin{multline*}
[K_{{\bf c}_2}{\bf k}_2,K_{{\bf c}_1}{\bf k}_1]-
[\mathcal{M}_S^*(K_{{\bf c}_2}{\bf k}_2),
\mathcal{M}_S^*(K_{{\bf c}_1}{\bf k}_1)]\\
=[K_{{\bf c}_2}({\bf c}_1){\bf k}_2, {\bf k}_1]-
\sum_{n=0}^{\infty}[{\overline{S({\bf
c}_2)}}^{(n)}{\overline{\bf c}_2}^{[n]}{\bf k}_2,{\overline{S({\bf
c}_1)}}^{(n)}{\overline{\bf c}_1}^{[n]}{\bf k}_1]\\
=\sum_{n=0}^{\infty}[{{\bf c}_1}^{[n]}\left(1-S({\bf
c}_1)^{(n)}{\overline{S({\bf
c}_2)}}^{(n)}\right){\overline{\bf c}_2}^{[n]}{\bf k}_2,{\bf k}_1]=
[\mathcal{K}_S({\bf c_1},{\bf c_2}){\bf k}_2,{\bf k}_1],
\end{multline*}
where \(\mathcal{K}_S\) is given by \eqref{defks}.
Since elements of the form
\[
F=\sum_{\ell=0}^m K_{{\bf c}_\ell} {\bf k}_\ell\]
are dense in ${\bf H}_2(\mathcal T)$,
\(S\) is a Schur multiplier if and only if
for any such element \(F\)
it holds that
\[[F,F]-[\mathcal{M}_S^*F,\mathcal{M}_S^*F]\ge 0,\]
which, in view of the computations above, is equivalent to the
positivity condition \eqref{poscon}.

\end{proof}

Below, following the analysis of  \cite[p. 86--90]{MR93b:47027},
we give an example of a Schur multiplier.
Let \({\bf a}\in {\mathbb D}({\mathcal T})\).
Then  \(K_{\bf a}({\bf a})\geq 0\) (in the sense of \(\mathbb K\))
and, since
\[K_{\bf a}({\bf a}) =1+\overline{\bf
a}{\bf a}
K_{\bf a}({\bf a})^{(1)},\]
  it  is also invertible. Hence the element
\[L_{\bf a}=K_{\bf a}({\bf a})^{(1)}K_{\bf a}({\bf a})^{-1}\]
 is also positive and
invertible. Since
\[ L_{\bf a}^{[k]}=K_{\bf a}({\bf a})^{(k)}K_{\bf a}({\bf a})^{-1},\]
 we have \(\rho(L_{\bf a})\leq 1\).

\begin{definition}
 The operator
\begin{equation*}
B_{\bf a}= ({\g}-{\bf a})(1-L_{\bf a}{\overline{\bf a}}\
{\g})^{-1}\sqrt{L_{\bf a}}\in{\bf U}(\mathcal T)
\end{equation*}
is called the  Blaschke factor, corresponding to \({\bf a}.\)
\end{definition}

\begin{proposition}\label{blaprop}
The multiplication operator ${\mathcal M}_{B_{\bf a}}$
is an isometry in \({\mathbf H}_2(\mathcal T).\)
\end{proposition}

\begin{proof}
First of all, we note that, according to \eqref{inprod},
\(\mathcal{M}_{\g}\) is an isometry in \({\mathbf H}_2(\mathcal T),\)
and that
\[L_{\bf a}^{(m)}=L_{{\bf a}^{(m)}},\ B_{\bf a}\g^m=\g^m B_{{\bf a}^{(m)}}.\]
Hence
it is enough to show that
for any \(m\in\mathbb{Z}^+\) and \({\bf p},{\bf q}\in\mathbb{K}_2\)
 \[[B_{\bf
a}\g^m{\bf p},B_{\bf a}{\bf q}]=\delta_{m,0}[{\bf p},{\bf q}].\]
To check this, we rewrite \(B_{\bf a}\) in the form
\[B_{\bf a}=\left(\g K_{{\bf a}^{(1)}}{{K_{{\bf a}}({\bf
a})}^{(1)}}^{-1}-{\bf a}\right)\sqrt{L_{\bf a}}.\]
 Then for \(m>0\) we have
\begin{multline*}
[B_{\bf a}{\g}^m{\bf p},B_{\bf a}{\bf q}]=[B_{{\bf
a}^{(1)}}{\g}^{m-1}{\bf p},K_{{\bf a}^{(1)}}{{K_{{\bf a}}({\bf
a})}^{(1)}}^{-1}\sqrt{L_{\bf a}}{\bf q}]
\\=[ {B_{{\bf a}^{(1)}}( {\bf a}^{(1)})}^{(m-1)} {{\bf a}^{(1)}}^{[m-1]}{\bf p},
 {{K_{{\bf a}}({\bf a})}^{(1)}}^{-1}\sqrt{L_{\bf a}}{\bf q}]= 0.
 \end{multline*} Analogously,
\begin{multline*}[B_{\bf a}{\bf p},B_{\bf a}{\bf q}]=
[K_{{\bf a}^{(1)}}{{K_{{\bf a}}({\bf a})}^{(1)}}^{-1}\sqrt{L_{\bf
a}}{\bf p},K_{{\bf a}^{(1)}}{{K_{{\bf a}}({\bf
a})}^{(1)}}^{-1}\sqrt{L_{\bf a}}{\bf q}]+[\overline{\bf a}{\bf
a}L_{\bf a}{\bf p},{\bf q}]
\\=[(L_{\bf a}{{K_{{\bf
a}}({\bf a})}^{(1)}}^{-1}+1-{{K_{{\bf a}}({\bf a})}}^{-1}){\bf
p},{\bf q}]=[{\bf p},{\bf q}].
\end{multline*}
\end{proof}

As a corollary of Lemma \ref{proppoint} and Proposition
\ref{blaprop} we obtain that an element $F\in{\mathbf H}_2({\mathcal T})$
vanishes at the point ${\bf c}\in\mathbb D({\mathcal T})$ if and only if it
can be written as $F=B_{\bf c}G$ where $G\in{\mathbf H}_2(
{\mathcal T})$  is such
that $[G,G]=[F,F]$. More generally, one can consider the following
homogeneous interpolation problem:\\

{\sl Given ${\bf c}_1,\ldots,{\bf c}_N\in{\mathbb D}({\mathcal
T})$ find all $F\in{\bf H}_2({\mathcal T})$ such that}
\begin{equation}\label{interpol}
F({\bf c}_j)=0,\quad j=1,\ldots, N.\end{equation}

We assume that it is possible to recursively define
invertible \({\bf k}_j\)  by
\[
{\bf k}_1=1,\ {\bf k}_{j+1}=\bigm(B_{{\bf k}_1^{(1)}{\bf
k}_1^{-1}{\bf c}_1} B_{{\bf k}_2^{(1)}{\bf k}_2^{-1}{\bf
c}_2}\cdots B_{{\bf k}_j^{(1)}{\bf k}_j^{-1}{\bf c}_j}\bigm)({\bf
c}_{j+1}).
\]
Then \(F\) is a solution of the interpolation problem
\eqref{interpol}
if and only if
\begin{equation}\label{intersol}
F=\bigm(B_{{\bf k}_1^{(1)}{\bf k}_1^{-1}{\bf c}_1}
B_{{\bf k}_2^{(1)}{\bf k}_2^{-1}{\bf c}_2}\cdots B_{{\bf
k}_N^{(1)}{\bf k}_N^{-1}{\bf c}_N} \bigm)G,\end{equation} where
\(G\in {\bf H}_2({\mathcal T})\) satisfies $[G,G]=[F,F]$.
Indeed, assume  that
\[ F=G_0=\ldots=B_{{\bf k}_1^{(1)}{\bf k}_1^{-1}{\bf c}_1}
B_{{\bf k}_2^{(1)}{\bf k}_2^{-1}{\bf c}_2}\cdots B_{{\bf
k}_n^{(1)}{\bf k}_n^{-1}{\bf c}_n}G_n,\]
where \(G_n\in {\bf H}_2({\mathcal T}).\)
Then
\begin{equation*}
0=f({\bf c}_{n+1})=({\bf k}_{n+1}G_n)({\bf c}_{n+1})= G_n({\bf
k}_{n+1}^{(1)}{\bf k}_{n+1}^{-1}{\bf c}_{n+1}),\end{equation*}
and hence
\begin{equation*} G_n=B_{{\bf k}_{n+1}^{(1)}{\bf k}_{n+1}^{-1}{\bf c}_{n+1}}
G_{n+1},
\end{equation*}  where
 \(G_{n+1}\in {\bf H}_2({\mathcal T}),\) and \eqref{intersol}
 follows by induction.

\bibliographystyle{plain}

\def\cprime{$'$} \def\cprime{$'$} \def\cprime{$'$}

\bigskip

\end{document}